\documentclass{amsart}
\usepackage[usenames,dvipsnames,svgnames,table]{xcolor}
\usepackage{amsfonts, amssymb}
\usepackage{enumerate}
\usepackage{graphicx}
\usepackage{amsmath}
\usepackage{enumitem}

\newtheorem{theorem}{Theorem}

\newtheorem{cor}[theorem]{Corollary}

\newtheorem{prop}[theorem]{Proposition}
\newtheorem{lem}[theorem]{Lemma}

\theoremstyle{definition}
\newtheorem{defn}[theorem]{Definition}

\newtheorem{remark}[theorem]{Remark}

\newcommand{\annr}{\operatorname{ann}_r}
\newcommand{\annl}{\operatorname{ann}_\ell}

\begin{document}
		\title{Characterizations in terms of Generalized Inverses}
	\author{Feroz Siddique}
	\address{Department of Mathematics, University of Wisconsin Eau Claire, WI-54701, USA}
	\email{siddiquf@uwec.edu}

	\begin{abstract}
 Using simply outer inverse and their differences, we characterize many known classes of rings including von Neumann regular, unit-regular, exchange, clean, and  strongly clean rings.
	\end{abstract}
	
	\maketitle
	
\section{Introduction}
All rings considered are associative rings with identity, and all modules are left modules unless otherwise specified. For a ring $R$, we denote by $U(R)$ the group of units of $R$, by $J(R)$ its Jacobson radical, by $N(R)$ the set of nilpotent elements of $R$, and by $C(R)$ the set of central elements of $R$. For idempotents, the notation $\leq$ denotes the natural partial order, \emph{i.e.} $e\leq f$ if and only if $ef=fe=e$.\\

The study of generalized inverses occupies a unique position at the intersection of ring theory, linear algebra, operator theory, and semigroup theory with significant applications in numerical analysis, statistics, optimization, and systems theory. Although subtle forms of generalized inverses had already appeared in the context of operator algebras, the original algebraic concepts of inner and outer inverses can be traced to E.H. Moore’s seminal work on generalized reciprocals for singular matrices \cite{Moore} in 1920. These ideas were rediscovered independently by Penrose in 1955, whose celebrated characterization of generalized inverses for matrices\footnote{As pointed out by R. Hartwig \cite{hartwig3}, the 1-2-3-4 conditions of the Moore-Penrose inverse $a^{\dagger}$ for an element $a$, was parallely shown to exist generally  by Kaplansky  \cite{kaplansky1} in 1955 for a regular ring $R$  with any proper involution $*$, such as one defined as  the conjugate transpose on a matrix ring over a complex field, and we therefore call it  the Moore-Kaplansky-Penrose Inverse.} \cite{penrose1} led to their systematic development and widespread application. Generalized inverses have since become indispensable in the study of singular linear systems, operator equations, and the algebraic structure of non-invertible elements.\\

We will first briefly recall some basic notions and  definitions of specific types of generalized inverses, namely  the sets of inner and outer inverses. An element $a \in R$ is called \textit{von Neumann regular} if the left (right) principal ideal generated by $a,$  is also generated by some idempotent $e \in R$. Equivalently, there exists  $b \in R$ such that \[a = aba \hspace{10mm} \cdots \hspace{3mm} (i)\] 
\noindent Here, $b$ is called an $1$\emph{-inverse} or  \emph{inner} of $a$.  In general,  a von Neumann regular element may have more than one inner inverse. The set of all inner inverses of $a$ is denoted by $I(a). $   
\noindent An element $a$ is called \emph{anti-regular} if there exists $b \in R$ such that \[b= bab \hspace{10mm} \cdots \hspace{3mm} (ii)\] Here $b$ is also called an \emph{outer-inverse} or \emph{pseudo-inverse} or $2$\emph{-inverse} of $a$. We will use pre-established notations to denote the set of all outer inverses of an element $a$ as $T(a)$ and the set of reflexive inverses as $\operatorname{Ref}(a)$.

\noindent It is remarkable that the two fundamental equations $(i)$ and $(ii)$ of generalized inverse theory did not originally arise as algebraic identities. Rather, they emerged as consequences of the existence of Moore’s \emph{general reciprocal} via projection identities. We will state the following theorem that explicitly denotes the existence of the outer and inner inverse, in its dense symbolic  form and is certainly not for the faint-hearted. 
%First we refer to A. Ben-Israel \cite{ben-israel} that restates Moore's 29.3 theorem of the existence of general reciprocal in modern notation.  

%\begin{theorem}[\cite{Moore} (29.3)]
	%
%	\[
%	\mathfrak{A}^{C}\mathfrak{P}^{1 \phantom{b} II}\mathfrak{P}^{2 \phantom{b} II}\,\kappa^{12}\cdot)\cdot	\]
%	\[
%	\exists |\,\lambda^{21
%	\text{ type }	
%	\mathfrak{M}^{2}_{\kappa^{*}}
%	\bar{\mathfrak{M}^{1}_{\kappa}};	}
%	\]
	
%	\[	
%	S^{2}\kappa^{12}\lambda^{21}	
%	=	
%	\delta^{11}_{\mathfrak{M}^{1}_{\kappa}},	
%	\qquad	
%	S^{1}\lambda^{21}\kappa^{12}	
%	=	
%	\delta^{22}_{\mathfrak{M}^{2}_{\kappa^{*}}}.	
%	\]
%\end{theorem}

\begin{theorem}[Moore (29.4)]
	\[	
	\mathfrak{A}^{C}\mathfrak{P}^{11}\mathfrak{P}^{22}\cdot	
	\lambda^{21}\!\cdot\kappa^{12}\cdot)\cdot	
	\]
	
	\[
	\begin{aligned}
		(1)\;&
		S^{1}S^{2}\lambda^{21}\kappa^{12}\lambda^{21}		
		=\lambda^{21},
		\qquad		
		S^{2}S^{1}\kappa^{12}\lambda^{21}\kappa^{12}		
		=\kappa^{12},\\[1ex]		
		(2)\;&		
		r_{\lambda}=r_{\kappa},\\[1ex]		
		(3)\;&		
		\mathfrak{M}^{1}_{\lambda^{*}}		
		=		
		\mathfrak{M}^{1}_{\kappa},		
		\qquad		
		\mathfrak{M}^{2}_{\lambda}		
		=
		\mathfrak{M}^{2}_{\kappa^{*}},\\[1ex]		
		(4)\;&		
		\lambda^{*12}
		=
		\operatorname{rec}\,\kappa^{*21},
		\qquad
		\kappa^{12}
		=
		\operatorname{rec}\,\lambda^{21}.
	\end{aligned}
	\]
	
\end{theorem}
In modern notation,  $\kappa^{12}$ is a finite matrix over a type C number system and $\lambda^{21}$ denotes  its unique  reciprocal(existence proved in Theorem 29.3 \cite{Moore}).  The symbols $r_\kappa$ and $r_\lambda$ denote the respective ranks, $(\cdot)^*$ denotes the adjoint transformation, and $\mathfrak M^i$ denotes Moore’s fundamental subspaces associated with the transformation.  Then:
\begin{enumerate}
	\item  The original transformation and its reciprocal satisfy the reciprocal identities
	$\kappa^{12}\lambda^{21}\kappa^{12}=\kappa^{12},
	\qquad
	\lambda^{21}\kappa^{12}\lambda^{21}=\lambda^{21},$\\
	where the operators $S^1S^2$ and $S^2S^1$  indicate  the transformations on the appropriate spaces.
	\textbf{In modern terminology, $\lambda^{21}$ is simultaneously an inner inverse and an outer inverse of $\kappa^{12}$;} %equivalently, it is a reflexive generalized inverse.
	\item The reciprocal has the same rank as the original matrix
	$r_\lambda=r_\kappa .$
	\item The corresponding fundamental subspaces coincide
	$	\mathfrak{M}^{1}_{\lambda^{*}}
	=
	\mathfrak M_\kappa^{\,1},
	\qquad
	\mathfrak M_{\lambda}^{\,2}
	=
	\mathfrak M_{\kappa^{*}}^{\,2}$
	that is, that the reciprocal preserves the appropriate range (equivalently, image) spaces associated with the transformation and its adjoint.
	\item Reciprocity is symmetric
	$\lambda^{*12}
	=
	\operatorname{rec}(\kappa^{*21}),
	\qquad
	\kappa^{12}
	=
	\operatorname{rec}(\lambda^{21})$
	that is, taking the reciprocal commutes with taking the adjoint, and the reciprocal of the reciprocal is the original transformation.
\end{enumerate}

For a reformulated account on E.H. Moore's pioneering work on generalized reciprocals for singular matrices, we refer to A. Ben-Israel expository article \cite{ben}.\\

%Let $\kappa^{12}$ be a finite matrix over a type C number system, and let $\lambda^{21}$ be its unique general reciprocal (existence proved in Theorem 29.3 \cite{Moore}). Then the following assertions hold:

% In such a case $b$ is also called an \emph{anti-regular} element. An element $x$ is called an \emph{inner inverse}  of an element $a$ if $axa = x$. 

\section{Outer and Inner Inverses}

\bigskip

The theory of generalized inverses plays a fundamental role in the study of noninvertible elements in rings. Among these, inner inverses and outer inverses occupy a central position. The existence of an inner inverse characterizes von Neumann regular elements, making generalized inverses an indispensable tool in the study of regularity. Beyond regular rings, outer inverses have developed into an important subject in their own right, with deep connections to matrix theory, operator theory, semigroup theory, and the structure theory of rings. Consequently, many questions concerning regularity, unit-regularity, exchange properties, and one-sided inverses may be reformulated in terms of the existence and structure of outer inverses.
For foundational definitions of outer inverses, we refer the reader to M.P. Drazin’s work~\cite{drazin1}.\\

In this section we review the known descriptions of the sets of inner and outer inverses. While the collection of all inner inverses admits an affine parametrization, no analogous linear expression is known for the collection of outer inverses in an arbitrary ring due to its complex algebraic variety. %\cite{predrag}. . This contrast provides the starting point for the results developed in this paper.\\

We begin with the well-known parametrization of the set of all inner inverses. For any $ a\in R$ and $x \in I(a),$  the class of all inner inverses of   $a$ is given as follows, as noted in  Rao and Mitra \cite{rao_mitra}.
\begin{enumerate}
	\item   $I(a) = x + (1-xa)b + c(1-ax)$ 
	\item $I(a) = x + b - xabax$
\end{enumerate}
for arbitrary $b,c \in R$.\\

This description stands in sharp contrast to the theory of outer inverses, where even the question of existence requires a more delicate characterization. The following  results and notation by Jain et al. indicates when an outer inverses of an element exists.

\begin{lem}[\cite{bapat}]
	Given an element $a$ in a semigroup $S$, the following statements are equivalent.
	\begin{enumerate}
		\item[(i)] $a$ has an outer inverse in $S$.
		\item[(ii)] There exists an idempotent element $e \in S$ such that $eS \subseteq aS$.
		\item[(iii)] There exists an idempotent element $f \in S$ such that $Sf \subseteq Sa$.
	\end{enumerate}
	%	In fact, $e$ in \emph{(ii)} is given by $ax$ for some $x \in \{a^{=}\}$ and $f$ in \emph{(iii)} is given by $ya$
In fact $e = ax = ya$	for some $x,y \in T(a)$.
\end{lem}
%
%\begin{proof}
%	Let $c$ be an arbitrary outer inverse of $a$. Now for $e=ac$, we observe that
%	$eS \subseteq aS$. This proves \emph{(i)}$\Rightarrow$\emph{(ii)}.

%	To prove the converse, let $e$ be an idempotent such that $eS \subseteq aS$.
%	Since $e$ is an idempotent, it follows that $e \in eS \subseteq aS$ and therefore
%	$e=az$ for some $z \in S$. By taking $x=ze$, observe that $xax=x$ and $ax=e$
%	hold by direct substitution. So, we have \emph{(ii)}$\Rightarrow$\emph{(i)}.

%	The proof of \emph{(i)}$\Leftrightarrow$\emph{(iii)} is similar to that of
%	\emph{(i)}$\Leftrightarrow$\emph{(ii)}.
%\end{proof}
%\bigskip 

\noindent Specialising this to the “set description” gives an explicit union:
\[
T(a)
=\bigcup_{\substack{e^2=e\\ e\in aR}}
\ \{\, ze \;:\; z\in R,\ az=e \,\}.
\]
By the left–right symmetric clause in the same characterisation, one also has
\[
T(a)
=\bigcup_{\substack{f^2=f\\ f\in Ra}}
\ \{\, f w \;:\; w\in R,\ wa=f \,\}.
\]
\bigskip 

%However, unlike the affine subspace of the set of all  inner inverses, an explicit linear expression  that characterizes outer inverses of an element in an arbitrary ring $R$, the form of a known outer inverse is missing in literature, due to its complex algebraic variety. %\cite{predrag}. 
The problem of describing all outer inverses of a fixed element has been studied in several important classes of rings. In general, no complete parametrization is known for arbitrary rings due to its complex variety, but explicit descriptions have been obtained under additional structural assumptions. These descriptions reveal how the geometry of idempotents and the regularity properties of the underlying ring determine the form of all solutions to the equation

\[
xax=x.
\]

\noindent Hartwig  obtained an explicit parametrization in the setting of anti-regular rings.
\begin{theorem}{\cite{hartwig3}}
	In an anti-regular ring, the general solution of $xax = x$ is given by 
	\[x = f(eaf)^{\hat{} } \phantom{l}e\]
	where $e$ and $f$ are arbitary idempotents and $(.)^{\hat{} } $ denotes any outer inverse.
\end{theorem}

\noindent The presence of inner inverses in regular rings allows one to derive a considerably more flexible description of all outer inverses. Rao and Mitra \cite{rao_mitra} proved the following characterization.
%On the other hand, when the element is regular, the additional existence of inner inverses yields a considerably richer structure, and the complete description of the set of outer inverses is given by Rao and Mitra in \cite{rao_mitra}, stated below. 
\begin{theorem}{ \cite{rao_mitra}}
	In a regular ring, the general solution of $xax = x$ is given by 
	\[x = b(cab)^{\hat{} } \phantom{l}b\]
	where $b$ and $c$ are arbitary elements and $(.)^{\hat{} } $ denotes any outer inverse.
\end{theorem}

Although Rao and Mitra require the ambient ring to be regular, the regularity assumption can in fact be localized to the element itself. Specifically, Patricio and Hartwig \cite{hartwig}  explicitly parametrized  the class of outer inverses for idempotents and group invertible elements. 

\begin{theorem}[ \cite{hartwig} ,Theorem 2.1]\label{thm}
	If $e$ is idempotent, then
	\[
	T(e) = v(x)\, g\, w(y),
	\]
	where 	$g\le e$, and 
	\[
	v(x)=1+(1-e)xe,
	\qquad
	w(y)=1+ey(1-e),
	\]
	for $x$ and $y$ are arbitrary.
\end{theorem}

\begin{cor}[\cite{hartwig}, Corollary 2.4]
	If $a$ has a group inverse $a^{\#}$, then
	\[
	T(a)=\bigl[1+(1-aa^{\#})x\,aa^{\#}\bigr]\;g\;\bigl[1+aa^{\#}y(1-aa^{\#})\bigr]\,(a^{\#}+1-aa^{\#}),
	\]
	where $g\le aa^{\#}$.
\end{cor}

Extending it to the case when $a$ is von Neumann regular, Lee and Lin \cite{leeandlin} provides a characterization for all outer inverses and reflexive inverses of $a$, using Theorem \ref{thm} above. 
%  Lee proved that whenever $a$ is a von Neumann regular element of an arbitrary ring with unity, the sets of all outer inverses and reflexive inverses of $a$ admit explicit parametrizations.

\begin{lem}[\cite{leeandlin}]
	Let $R$ be a ring and let $a \in \text{Reg}(R)$, $a^{-} \in I(a)$ and $e := aa^{-}$. Then the following hold:
	\begin{enumerate}[label = (\roman*)]
		\item $T(a)=\Bigl\{\bigl(a^{-}+(1-a^{-}a)r\bigr)\,v(x)\,g\,w(y)\ \bigm|\ r,x,y\in R,\ g=g^{2}\in eRe\Bigr\}$
		
		\item $\text{Ref }(a)=\Bigl\{\bigl(a^{-}+(1-a^{-}a)r\bigr)\,v(x)\,e\,w(y)\ \bigm|\ r,x,y\in R\Bigr\}$
		
	\end{enumerate}
	where
	\[
	v(x)=1+(1-e)xe
	\quad\text{and}\quad
	w(y)=1+ey(1-e)
	\qquad (x,y\in R).
	\]
\end{lem}

\noindent On careful observation, we see that the factor $v(x)$ appearing in Lee's parametrization can be refined further as we see below. 
%\[ \text{Ref }(a)=\Bigl\{\bigl(a^{-}+(1-a^{-}a)r\bigr)\,\,e\,w(y)\ \bigm|\ r,x,y\in R\Bigr\} \] 
%\[T(a)=\Bigl\{\bigl(a^{-}+(1-a^{-}a)r\bigr)\,\,g\,w(y)\ \bigm|\ r,x,y\in R,\ g=g^{2}\in eRe\Bigr\}\]
%where $w(y)=1+ey(1-e)  \text{ for }(x \in R).$

\begin{lem}{\label{refinement}}
	Let $R$ be a ring and let $a \in \text{Reg}(R)$, $z \in I(a)$ and $e = az$. Then the following hold:
	\begin{enumerate}[label = (\roman*)]
		\item $T(a)=\Bigl\{\bigl(z +(1-za)r\bigr)\,\,g\,w(y)\ \bigm|\ r\in R,\ g=g^{2}\in eRe\Bigr\}$
		
		\item $\text{Ref }(a)=\Bigl\{\bigl(z +(1-za)r\bigr)\,\,e\,w(y)\ \bigm|\ r \in R\Bigr\}$
		
	\end{enumerate}
	Here $w(y)=1+ey(1-e) \qquad (x,y\in R).$ 
	
\end{lem}
\begin{proof}
$(i)$ \quad 	By Lemma 7, we have  \[T(a) = \Bigl\{\bigl(z +(1-za)r\bigr)\,v(x)\,g\,w(y)\ \bigm|\ r\in R,\ g=g^{2}\in eRe\Bigr\}\] Fixing the factor $v(x)$ to be the identity, we  have \[\Bigl\{\bigl(z +(1-za)r\bigr)\,\,g\,w(y)   \ \bigm|\ r\in R,\ g=g^{2}\in eRe\Bigr\} \subseteq T(a)\]

\noindent 	Conversely, let $x\in T(a)$.  By Lee's characterization, there exist 	$r,s,y\in R$ and an idempotent $g\in eRe$ such that
	\[	
x=\bigl(z+(1-za)r\bigr)v(s)\,g \,w(y),	
	\]
	where $v(s)=1+(1-az)s(e)$ and $	w(y)=1+azy(1-az).$ \\
	Then we can write 
	\begin{align*}
		x  &= \bigl(z+(1-za)r\bigr)\, \underbrace{(1+(1-az)s(e))}_{v(s)} \,g \,w(y) \\
		&= \bigl(z + (1-za)r + (z-zaz)se + (1-za)r(1-az)se\bigr)\,g \,w(y) \\
		&= \bigl(z + (1-za)(r +zse + r(1-az)se)\bigr)\,g \,w(y) \\
		&= \bigl(z + (1-za)r' \bigr)\,g \,w(y) \text{ for some } r' \in R\
	\end{align*}
\end{proof}
%Let us consider the set $S = \Bigl\{\bigl(a^{-}+(1-a^{-}a)r\bigr)\,\,g\,w(y)\ \bigm|\ r,x,y\in R,\ g=g^{2}\in eRe\Bigr\}$ where $w(y)=1+ey(1-e)  \text{ for }(x \in R).$
The refinement developed above is quite useful as we shall see in subsequent results. One of its uses is that it  allows us to relate arbitrary outer inverses of a regular element to its reflexive generalized inverses through the natural partial order.  
%The significance of the  refinement is that it is more feasable to obtain the following relationship between a reflexive generalized inverse of a regular element to its outer inverse. The following theorem was conjectured in  Jain et.al \cite{jain} and shown true for any von Neumann regular element in an associative ring. \\
To state these results, we first recall the minus (or natural) partial order introduced independently by Hartwig \cite{hartwig1} and Nambooripad \cite{namboor}.
\begin{defn}
	Let $a,b\in R$ be von Neumann regular elements. Then
	\[
	a\leq^- b
	\]
	if there exists $x\in I(a)$ such that
	\[
	ax=bx
	\qquad\text{and}\qquad
	xa=xb.
	\]
\end{defn}

\noindent The minus order $a \leq^{-} b$ admits several  equivalent descriptions that are used in this paper. 
\begin{theorem}{\label{partialorder}}
	Let $a,b\in R$ be von Neumann regular. The following are equivalent.
	\begin{enumerate}[label=(\roman*)]
		\item $a\leq^- b$.
		\item \cite{gomes} There exists $x_r\in Ref(a)$ such that
		\[
		a=ax_ra=bx_ra=ax_rb.
		\]
	\item \cite{hartwig1} There exists $x_r\in Ref(a)$ such that
		\[
		ax_r=bx_r
		\qquad\text{and}\qquad
		x_ra=x_rb.
		\]
		\item \cite{namboor} There exist idempotents $e,f\in R$ such that
		\[
		a=eb=bf.
		\]
		\item \cite{hartwig1}  There exists $a^{+} \in \operatorname{Ref}(a)$ and $s \in R$ such that 
		\[b = a+ (1-aa^{+})s(1-a^{+}a)\]
	\end{enumerate}
	\end{theorem}

The following connection was conjectured by Jain \emph{et al.} \cite{bapat}, who established it for von Neumann regular elements in associative rings. \\

\begin{theorem}{\cite{bapat}}
	Let $a$ be a (von Neumann) regular element in an associative ring $R$ and  $x$ be an outer inverse of $a$. Then there exists a reflexive  inverse $y$ of $a$ such that $x \leq^{-} y$. \\
	%	Furthermore, $x_r -x$ is von Neumann regular and $(x_r-x )\leq^{-} x_r$.
\end{theorem}

Using the refined Lemma \ref{refinement}, the result  below  recovers this theorem by a direct argument, and also demonstrates that the reflexive inverse dominating a given outer inverse is far from unique and satisfies several additional structural properties. \\ 
\begin{theorem}\label{outer0}
	Let $R$ be a ring, let $a$ be a von Neumann regular element of $R$, and let $x$ be an arbitrary outer inverse of $a$. Then the following hold.
	\begin{enumerate}[label=(\roman*)]
		\item There exists a family of reflexive inverses $\{x_{r,k}\}$ of $a$ such that
		\[
		x\leq^- x_{r,k}
		\]
		for every $k$.
		\item Every member of the family,
		\[
		\{x_{r,k}-x\}
		\]
		is von Neumann regular with a common inner-inverse and
		\[
		x_{r,k}-x\leq^- x_{r,k}.
		\]
	\end{enumerate}
\end{theorem}

%\begin{theorem}\label{outer}
%	Let $R$ be a ring and  $a$ be a von Neumann regular element in $R$ and  $x$ is an arbitrary outer inverse of $a$. Then the following conditions hold 
%	\begin{enumerate}[label = (\roman*)]
	%	\item There exists a family of reflexive  inverses $\{x_{r,k}\}$ of $a$ such that $x \leq^{-} x_{r,k}$. 
	%	\item 	Each member of the family  $\{x_{r,k} -x\}$ is von Neumann regular and $(x_{r,k}-x )\leq^{-} x_{r,k}$.
%	\end{enumerate}
%\end{theorem}
\begin{proof}
	Let $a = aza$ for some  $z \in R$ and let $x$ be any  outer inverse of $a$. Using the refinement in Lemma \ref{refinement}, we can write 
	\[x = \left(z  + (1-za)r_{\scriptscriptstyle 0}\right) g v \hspace{2mm} \text{ for some fixed } r_{\scriptscriptstyle 0} \in R \text{ and } g^2=g \in (az)R(az)\]
	where $v = 1 + (az)t(1-az)$ for some fixed $t \in R$. \\
	
	Set $x_{s,k}= \left[z  + (1-za)(r_{\scriptscriptstyle 0} + s(az-g))\right] (az) v_{\scriptscriptstyle k}$ where $v_{\scriptscriptstyle k} = 1 + az(gt +(az-g)k)(1-az)$ for arbitrary $s,k \in R$. Clearly $v_{\scriptscriptstyle k} \in U(R).$\\
	Then $a(x_{s,k}) = (az)v_{\scriptscriptstyle k}$ and $(x_{s,k} )a = \left[z  + (1-za)(r_{\scriptscriptstyle 0} + s(az-g))\right] a$. Then it can be easily seen that  $(x_{s,k})a(x_{s,k})=x_{s,k}$ and $a(x_{s,k})a= aza=a$. \\
	
	\noindent We prove that $x \leq^{-} x_{s,k}$. \\
	
	\noindent Take the element $ y = ga$ where $g^2=g \in (az)R(az)$ in the outer inverse $x = \left(z  + (1-za)r\right) g v_1 $.\\
	Since $vg = g = v_{\scriptscriptstyle k}g$ and $g^2=g$, we have 
	\begin{align*}
		x(ga)x = xg(ax) &= (xg )(gv)\\
		&= ( \left(z  + (1-za)r_{\scriptscriptstyle 0}\right) (g vg ) (gv) \\
		&=  \left(z  + (1-za)r_{\scriptscriptstyle 0}\right) g v  = x
	\end{align*}
	
	\noindent  Thus $ga$ is an inner inverse of $x$. Furthermore we have, 
	\begin{align*}
		(x_{s,k})(ga)  &= \left[z  + (1-za)(r_{\scriptscriptstyle 0} + s(az-g))\right] (az) v_k (ga)\\
		&=  \left[ z  + (1-za)(r_{\scriptscriptstyle 0} + s(az-g))\right] (ga)\hspace{4mm} \cdots \hspace{4mm} \because (az) g = g \\
		&=  \left[z  + (1-za)(r_{\scriptscriptstyle 0})\right] (ga) \\
		&= x(ga)
	\end{align*}
	and \hspace{2mm}	$(ga)(x_{s,k})  = g(a(x_{s,k}))=  g (az)v_{\scriptscriptstyle k} = g(v_{\scriptscriptstyle k})$ \hspace{1mm} where 
	%	&=  g (az) v_1 \hspace{4mm} \cdots \hspace{4mm} \because g (az) = g \\
	%	&=  g(ax_r ) = (ga)x_r
	\begin{align*}
		g  (v_{\scriptscriptstyle k}) &= g(1 + az(gt +(az-g)k)(1-az)) \\
		&= g + g(gt +(az-g)k)(1-az)) \\
		&= g + gt (1-az) \\
		&= g (az) + g(az)^2t (1-az) \hspace{4mm} \cdots  \hspace{4mm} \because g(az) = g = g(az)^2\\
		&= g(az)(1 + (az)t(1-az)) = g(1 + (az)t(1-az)) =  gv.
	\end{align*}
	Therefore, $(ga)(x_{s,k})  = gv = g(ax) = (ga)x$ \\
	Thus,  $x \leq^{-} x_{s,k}$ \\
	\noindent $(ii)$ We now prove that each member of the family  $\{x_{s,k} -x\}$ is von Neumann regular and $(x_{s,k}-x )\leq^{-} x_{s,k}$.\\
	First, it can be seen that $a(x_{s,k} - x) = (az)v_{\scriptscriptstyle k}  - gv  $ and 
	\begin{align*}
		(x_{s,k} - x) (az - g) &= x_{s,k} (az) - (x_{s,k}) g - x(az) + xg \\
		&= \left[z  + (1-za)(r_{\scriptscriptstyle 0} + s(az-g))\right] (az)  - \left[z  + (1-za)(r_{\scriptscriptstyle 0} + s(az-g))\right] g  \\ & \hspace{4mm} - \left[z  + (1-za)(r_{\scriptscriptstyle 0})\right] g + xg \\
		&= (x_{s,k})v_{\scriptscriptstyle k}^{-1} - \left[z  + (1-za)(r_{\scriptscriptstyle 0})\right] g  \\
		\implies (x_{s,k} - x) (az - g) v_{\scriptscriptstyle k} &=  (x_{s,k}) - \left[z  + (1-za)(r_{\scriptscriptstyle 0})\right] gv_{\scriptscriptstyle k}
	\end{align*}
	\begin{align*}
		\text{Therefore,  } (x_{s,k} - x) (az - g)a (x_{s,k} - x) &= (x_{s,k} - x) (az - g) a (x_{s,k} - x) )\\
		&= (x_{s,k} - x) (az - g)((az)v_{\scriptscriptstyle k}  - gv) \\
		&= (x_{s,k} - x) ((az)v_{\scriptscriptstyle k}  - gv - gv_{\scriptscriptstyle k} + gv)\\
		&= (x_{s,k} - x) (az - g) v_{\scriptscriptstyle k} \\
		&= (x_{s,k})- \left[z  + (1-za)(r_{\scriptscriptstyle 0})\right] (g v_{\scriptscriptstyle k}) \\
		&= (x_{s,k})- \left[z  + (1-za)(r_{\scriptscriptstyle 0})\right] (g v) \\
		&= (x_{s,k} - x)
	\end{align*}
	Thus, each member of the family  $\{x_{s,k} -x\}$ is von Neumann regular with a common inner inverse $a-ga$.
	\begin{align*}
		\text{Furthermore,  }  (az - g)a (x_{s,k} - x) &= (az - g)((az)v_{\scriptscriptstyle k}  - gv) \\ 
		&= (az - g) v_{\scriptscriptstyle k} \\ 
		&= (az - g) (az)v_{\scriptscriptstyle k} \\ 
		&= (az - g) a(x_{s,k}) \hspace{50mm}
	\end{align*}
	Now to show that $(x_{s,k} - x)  (az - g)a  = (x_{s,k} )(az - g)a  $ we first observe that since $(a-ga) = (az-g)(a-ga)$ we have 
	\begin{align*}
		(x_{s,k} )(az-g)a &=  \left[z  + (1-za)(r_{\scriptscriptstyle 0} + s(az-g)\right] (az-g) a  \\ 
		&= \left[z  + (1-za)(r_{\scriptscriptstyle 0} + s)\right]  (az-g) (a-ga)  \\ 
		&= \left[z  + (1-za)r_{\scriptscriptstyle 0} \right](a-ga)   + \left[(1-za)s\right](az-g)(a-ga)   \\
		&=  \left[z  + (1-za)r_{\scriptscriptstyle 0} \right]a -  \left[z  + (1-za)r_{\scriptscriptstyle 0} \right](ga) + \\ 
		& \hspace{4mm} \left[(1-za)s\right](az - g)(a) - \underbrace{\left[(1-za)s\right](az - g)(ga) }_{0} \\
		&=  \left[z  + (1-za)(r_{\scriptscriptstyle 0} + s(az-g)\right] a  -   \left[z  + (1-za)r_{\scriptscriptstyle 0} \right](ga) 
	\end{align*}
	\begin{align*}
		\text{ where, }  (x_{s,k} - x)(az-g)a &= (x_{s,k})v_{\scriptscriptstyle k}^{-1}a - \left[z  + (1-za)(r_{\scriptscriptstyle 0})\right]  ga \\ 
		&=  \left[ z  + (1-za)(r_{\scriptscriptstyle 0} + s(az-g)\right] (a) - \left[z  + (1-za)(r_{\scriptscriptstyle 0})\right]  ga \\ 
	\end{align*}
	Therefore, $(x_{s,k}-x )\leq^{-} x_{s,k}$ for each $s,k \in R$.
\end{proof}

As a consequence, we have the following important corollary that characterizes a unit-regular element in terms of its outer inverses. To do that we use an established characterization \cite{siddique1} of unit-regular elements $(i) \Longleftrightarrow (iii)$ below on inner inverses. %Th result is not new and is also proved in \cite{siddique1} and \cite{lam3}.

\begin{cor}{\label{outer}}
	Let $R$ be a ring and $a \in R$. The following are equivalent. 
	\begin{enumerate}[label = (\roman*)]
		\item $a$ is unit-regular
		\item  Every outer inverse  of $a$ is unit-regular. 
		\item $a$ has a unit-regular inner inverse.
	\end{enumerate}
\end{cor}
\begin{proof}
	$(iii) \implies (i)$ If $a$ has an inner inverse $y$ that is unit-regular,  then from \cite{lam} and  \cite{siddique1}, $a$ is unit regular. For the sake of completion, we provide an alternate proof below. \\
	Suppose $yuy = y$ for some unit $u \in U(R)$. Then $1 - uy(au^{-1})(1-uy) \in U(R)$ where 
	\begin{align*}
		1 - uyau^{-1}(1-uy) &=  1 - uya(u^{-1} - y)\\
		&= 1 - uya(1-yu)u^{-1}\\
		&= u\left (1 - ya(1-yu)\right )u^{-1}\\
		&= u\left (yu + (1-ya)(1-yu)\right )u^{-1}\\
		&= u\left (y+ (1-ya)(1-yu)u^{-1}\right )
	\end{align*}
	Therefore $y+ (1-ya)(1-yu)u^{-1}  = v \in U(R)$. Pre-multiplying both sides by $a$ gives us $ay = av$. Combined with $aya = a$ gives $a = ava$.\\
	$(i) \implies (iii)$ \emph{Straightforward}. \\
	
	\noindent 	$(i) \implies (ii)$ Let $a$ be unit-regular. Then by $(i) \Leftrightarrow (iii)$ every reflexive inverse $x_{\scriptscriptstyle r}$ is unit regular since it has $a$ has an inner inverse. \\
	Choose an arbitrary outer inverse $x \in \operatorname{T}(a)$. Then by Theorem \ref{outer0}, there exists a reflexive inverses $x_r$ of $a$ such that $x \leq^{-} x_r$ where $x_r$ is unit regular. We show that $x$ is unit-regular. \\
	Now that $x_r$ is unit-regular, we can write $x_r = e u$ for some idempotent $e^2 = e$ and $u \in U(R)$. Then $u^{-1}e$ is a reflexive inverse of $x_r$ since 
\[ x_r(u^{-1}e)x_r = eu(u^{-1}e)eu = eu = x_r \]

\[\text{and, } u^{-1}e(x_r)u^{-1}e = u^{-1}e(eu)u^{-1}e = u^{-1}e\]
Consequently, $x_r \leq^{-1} u$.\\
	By transitivity of the partial order, $x \leq^{-} x_r$ and $x_r \leq^{-} u \implies x \leq^{-} u $ for some unit $u \in U(R)$. Thus, there exists  $b \in  \operatorname{Ref}(x)$ such that $bx = bu \implies b = bxb = bub.$ Thus $x$ has an inner inverse $b$ which is unit regular, and consequently, $x$ is also unit-regular. %Subsequently, $x = xux$ for some unit $u$ in $U(R)$. 
	
	\noindent 	$(iii) \implies (i)$ If every outer inverse  of $a$ is unit-regular, then so is every reflexive inverse of $a$ and consequently, $a$ has an inner inverses $y$ that is unit-regular. 
\end{proof}

\noindent The result above strengthens the following proposition from   Hartwig \cite{hartwig1} from all   reflexive inverses to all outer inverses. 
\begin{prop}[\cite{hartwig1} Proposition 3 (iv)]
Let $R$ be a ring with unity and let $a$ be a von Neumann regular element of $R$. Then the following hold 
\begin{enumerate}
	\item $a$  is unit-regular  $\Longleftrightarrow a^{+}$ is unit regular for some and hence all reflexive inverses $a^{+}$ of $a$. 
\end{enumerate}
\end{prop}

We will now continue our journey and  provide additonal characterizations of unit-regularity in terms of its outer inverses, and in the process strengthen two additional conditions given in Hartwig \cite{hartwig1} that are stated below.
\begin{prop}[\cite{hartwig1} Proposition 3 (v),(vi)]
	Let $R$ be a ring with unity and let $a$ be a von Neumann regular element of $R$. Then the following hold 
	\begin{enumerate}[label = (\roman*)]
				\item $a$  is unit-regular  $\Longleftrightarrow$ there exists a reflexive inverses $x$ of $a$ and $t\in R$ such that 
		\[x +(1-xa)s(1-ax)) \in U(R)\]
				\item $a$  is unit-regular  $\Longleftrightarrow$ there exists a reflexive inverses $x$ of $a$ and $s\in R$ such that 
		\[a+(1-ax)s(1-xa)) \in U(R)\]
	\end{enumerate}
\end{prop}
\medskip

We will strengthen  the consequences and prove they are true for \textbf{all} outer inverses $x \in T(a)$ and all reflexive inverses $x \in Ref(a)$ respectively. Furthermore, we see that $s,t \in R$ can be taken as units. %Most importantly, we will soon see that the terms $(1-aa^{+})$  and $(1-a^{+}a)$  can be switched in the above decomposition, which is most peculiar and is not at all evident from direct calculation.
%But first we will give a direct ring-theoretic proof to following lemma of Ehrlich \cite{ehrlich} and Handelman's \cite{handelman} characterization of unit regular rings, and in that process weaken their assumption 

\begin{theorem}{\label{urr}}
Let $R$ be a ring with unity and let $a \in R$ be von Neumann regular. The the following are equivalent:
	\begin{enumerate}[label = (\roman*)]
		\item $a$ is unit-regular. 
	\item For every outer inverse $x \in T(a),$ there exists $s\in U(R)$ such that \[x+ (1-xa)s(1-ax) \in U(R) \] 
	\item For every reflexive inverse $x \in Ref(a)$, there exists $ t\in U(R)$ such that  \[a + (1-ax)t(1-xa)) \in U(R) \]
\end{enumerate}
\end{theorem}
\begin{proof}
	$(i) \implies (ii)$. Let $a$ be unit-regular. \\Choose an arbitrary outer inverse $x \in T(a)$. Then by Corollary \ref{outer}, $x$ is unit-regular and subsequently, there exists $v \in U(R)$ such that $x = xvx$. Label $e = ax, f = xa$. 
	\[\text{Set } \,\,c = (1-xa)v^{-1}(1-ax) \,\,\, \text{ and } \, d = v - vxv\]
Label	$p=xv$, $q=vx$. Then $p^2=p$, $q^2=q$, and  $d=(1-q)v=v(1-p)$. \\
In fact, we  have 
\[fp=(xax)v=xv=p \quad \text{ and } \quad  pf=(xvx)a=xa=f\]
\[eq=a(xvx)=ax=e \quad \text{ and } \quad qe=v(xax)=vx=q\]
Passing to complements and expanding similarly as above we also have  
\[(1-f)(1-p)=1-f \quad \text{ and } \quad  (1-p)(1-f)=1-p \] \[ (1-e)(1-q)=1-q \quad \text{ and } \quad  (1-q)(1-e)=1-e \]

\noindent Then using the above expansions we have
\[(1-e)d=(1-e)(1-q)v=(1-q)v=d \] and \[d(1-f)=v(1-p)(1-f)=v(1-p)=d \]
Subsequently, expanding $cd$ and $dc$ we get \\
\begin{minipage}{0.48\textwidth}
	\begin{align*}
		cd &= \left((1-f)v^{-1}(1-e) \right) \left(d\right) \\
		&=  \left((1-f)v^{-1} \right)\underbrace{(1-e)d}_{d} \\ 
		&=  \left((1-f)v^{-1} \right) v(1-p)\\
		&=(1-f)(1-p) \\
		&=1-f
	\end{align*}
\end{minipage}
\hfill
\begin{minipage}{0.48\textwidth}
	\begin{align*}
		dc &= (d)\left((1-f)v^{-1}(1-e) \right) \\
		&=  \underbrace{d(1-f)}_{d}v^{-1}(1-e)\\
		&=(1-q)v \left(v^{-1}(1-e)\right) \\
		&=(1-q)(1-e) \\
		&= 1-e\\
	\end{align*}
\end{minipage}

\noindent Then clearly, $(x+c) \in U(R)$ with inverse $axa + d,$  since  
\begin{align*}
		(axa+d)(x+c) &= (ax)(ax)   \, +  \, (axa)c  \, + \,  dx  \, + \, dc \\
		&= \underbrace{ax}_{e} + a\underbrace{(xa)c}_{0}  \, +  \, \underbrace{(v-vxv)x}_{0}  \, +  \, \underbrace{dc}_{1-e}=  1
\end{align*}
\begin{align*}
	(x+c)(axa + d) &= (xa)(xa)  \, +  \, c(axa) \, +  \, xd  \, +  \, cd \\
	&= \underbrace{xa}_{f} + \underbrace{c(ax)}_{0}a + \underbrace{x(v-vxv)}_{0} + \underbrace{cd}_{1-f}=  1
\end{align*}
Thus, $x+ (1-xa)s(1-ax) \in U(R) $ for some $s = v^{-1} \in U(R)$ where $xvx = x$.\\
%\begin{align*}
%	(axa+d)(x+c):&\quad (axa)x=a(xa\cdot x)=ax=e,\quad (axa)c=a\underbrace{(xa)(c)}_{0}=0,\\
%	&\quad dx=0,\quad dc=1-e;\\
%	(x+c)(axa+d):&\quad x(axa)=(xa)=f,\quad xd=\underbrace{x(1-e)}_{0}d=0,\\
%	&\quad c(axa)=c(ax)a=0,\quad cd=1-f,
%\end{align*}
% and $(axa+d)(x+c)= e+(1-e)=1=f+(1-f)=(x+c)(axa+d)$.\\
 
\noindent $(i) \implies (iii):$ Let $a$ be unit-regular. Choose an arbitrary  reflexive inverse $ x \in Ref(a)$ such that $axa = a$ and $xax = x$. Then $x$ is also unit regular  and there exists $v \in U(R)$ such that $x = xvx$. Then similarly as above, setting $e=ax, f=xa$ and $d = v-vxv$, we have $(axa +d)= a+d$ is invertible with inverse $x +c$. Since we also have  $(1-e)d = d$ and $d(1-f) = d,$ then together with $a = axa$ we get 
 \begin{align*}
 	(1-e)(a+d)(1-f) &= (1-e)a(1-f) \, + \, (1-e)d(1-f)\\
 	&= 0 \, + \,  (1-e)d(1-f) = d
 \end{align*}
Setting $t= (a+d)$ we get  $a + (1-ax)t(1-xa) \in U(R)$ for $t \in U(R)$. \\
$(ii) \implies (i):$ Let $x$ be any outer inverse of $a.$ Then by Theorem \ref{outer0}, there exists a reflexive inverse $y \in Ref(a)$ such that $x \leq^{-} y$. As $y$ is also an outer inverse of $a$,  there exists $s \in U(R)$ such that \[y + (1-ya)s(1-ay) \in U(R)\] 
whence  $a\left(y + (1-ya)s(1-ay) \right)a = aya = a.$ \\
$(iii) \implies (i):$ Choose an outer inverse $x$ of $a$. Then by Theorem \ref{outer0}  there exists a reflexive inverse $y \in Ref(a)$ such that $x \leq^{-} y$. By assumption, there exists $s \in U(R)$ such that \[a + (1-ay)s(1-ya) \in U(R)\] Then $y$ is unit-regular since $y(a + (1-ay)s(1-ya))y = y,$ and subsequently, $a$ has a unit-regular inner inverse $y$ and by Theorem \ref{outer}, $a$ is unit-regular.
\end{proof}

\begin{remark}
	We draw attention to the structure of the proof of $(i)\Rightarrow(ii)$ 
	above. Recall the classical criterion of unit regularity in Ehrlich \cite{ehrlich} (see also 
	Handelman \cite{handelman}):\\
	 \emph{A regular ring R is unit-regular if and only
	if for every two idempotents, $e \sim f \implies (1-e) \sim (1-f)$.}\\
The construction $c\in(1-f)R(1-e)$ and 
$d\in(1-e)R(1-f)$ satisfying $cd=1-f$ and $dc=1-e$  yields a purely element-wise, ring-theoretic proof of the classical characterization without module-theoretic machinery.
\end{remark}
%As a consequence we provide a characterization of stability range and strengthen Proposition 9 in \cite{siddique}.\\

\noindent Nicholson, in \cite{ASring}, introduced the following notion: 
\begin{defn}
An element $a \in R$ is said to be \emph{left annihilator-stable} (or simply, a \emph{left AS element}) if  whenever $Ra + \mathrm{ann}_l(b) = R$,  there exists a unit $u \in U(R)$ such that \[a - u \in \mathrm{ann}_l(b)\] 
A ring $R$ is called a \emph{left AS ring} if every element of $R$ is left AS.
\end{defn}

\noindent  He  further showed that for any element $a \in R$, the following implications hold:

\[
a \text{ is unit-regular } \;\implies\; a \text{ has Stable Range 1 } \;\implies\; a \text{ is AS}.\\
\]
\medskip 
Canfell gave the following characterization for left (right) AS rings.
\begin{theorem}[Canfell {\cite[Cor.~4.4]{canfell}}]\label{thm:canfell}
The following statements are equivalent for a ring $R$.	\begin{enumerate}[label=(\roman*)]
		\item If $Ra+ann_l(b)=R$ for $a,b\in R$, then there exists a unit $u\in U(R)$ such that $a-u\in ann_l(b)$.
		\item $R$ is left uniquely generated (left UG).
	\item If $Ra=Rb$ for $a,b\in R$, then $a=vb$ for some left unit $v\in R$.
		\end{enumerate}
\end{theorem}
\medskip 
We can  see several examples of AS rings, namely domains, local rings \cite{kaplansky}, unit-regular rings \cite{lam4}, in particular, rings with stable range one,  matrix rings over elementary divisor rings \cite{kaplansky}, $\begin{pmatrix}
	\mathbb{Z}_2 & M\\
	0 & \mathbb{Z}
\end{pmatrix}$
where  $M$ is a $(\mathbb{Z}_2,\mathbb{Z})$ -bimodule \cite{zhou},  etc. 
\noindent We wonder  for any element $a \in R,$ if the initial condition $Ra + \mathrm{ann}_l(b) = R, $ can be weakened $b \in R.$   If $1 = ra + y$ for some $r\in R$ and $y \in \operatorname{ann}_{l}(b)$, then it can is easily seen \[Ra + \mathrm{ann}_l(b) = R \iff b \in R(ab) \]In fact, if $a$ admits outer inverses then  for any $x \in T(a)$ and $b= xr \in xR$ we must have $x(ab) = xa(xr) = (xax)r= xr = b$. \\  Thus $\{b : Ra + \annl(b) = R\} \supseteq \bigcup_{x \in T(a)} xR$. \\

On the other hand, in addition to  $Ra + \mathrm{ann}_l(b) = R$ holding true for some $b \in R,$ if $ab$ is assumed to be von Neumann regular, then using the fact that 
$Ra + \mathrm{ann}_l(b) = R \iff b = r(ab)$ for some $r \in R,$ and $(ab)s(ab) = ab$ for some $s\in R$ we have 
\[(bs)a(bs) = (\underbrace{rab}_{b}s)a(bs)= r(ab)s(ab)s= rabs = bs \]
Consequently, $a$ admits an outer inverse $bs \in \operatorname{T}(a)$ such that $bs(ab) = (rab)s(ab) = rab = b \implies b\in xR. $ Thus $b \in \bigcup_{x \in T(a)} xR$.  \\
Infact, we have the following lemma.
\begin{lem}{\label{ASouter}}
	Let $R$ be an exchange ring. Then for any element $a \in R$, 
\[	\{b : Ra + \annl(b) = R\} \,\, = \,\,  \bigcup_{x \in T(a)} xR\]
\end{lem}
\begin{proof}
%The inclusion $\supseteq$ holds in any ring, as shown above.\\	
We prove that $ 	\{b : Ra + \annl(b) = R\} \subseteq \bigcup_{x \in T(a)} xR$.\\
Let
$Ra + \annl(b) = R \iff b = rab$ for some $r \in R$.  Label $a' = ra$ where $(1-a')b = 0$. Since $R$ is an exchange ring, by Nicholson's
characterization there is an idempotent $e \in Ra' \subseteq Ra$ with $1 - e \in R(1-a').$  Set idempotent $e = ra$ and  $1 - e = s(1-a')$. Then $(1-e)b = s(1-a')b = 0 \implies  eb = b $. \\
 Label $x = rar$. Then $xa  = e^2 = e$ and $xax = ex = x \implies x \in T(a)$ where
\[ x(ab) = (xa)b = eb = b   \]
Thus  $b \in xR$ for some $x \in T(a)$. \\
\end{proof}
\begin{remark}{\label{ASouter1}}
It is quite straightforward to see the right analogue 	\[\{b : aR + \operatorname{ann}_{r}(b) = R\} \,\, = \,\,  \bigcup_{x \in T(a)} Rx\]  holds true for exchange rings as well, using Nicholson's right ideal characterization and the fact that $\{b : aR + \operatorname{ann}_{r}(b) = R \}  \iff b \in (ba)R$
\end{remark}
\medskip

\noindent We can now see that the forward implication  above
\[
a \text{ is unit-regular } \;\implies\; a \text{ has Stable Range 1 } \;\implies\; a \text{ is AS}.\\
\]
 invites a  formal weakening. %If $a$ admits an outer inverse $x$, then $Ra+\mathrm{ann}_l(x)=R$ holds automatically. 
We weaken the left (right) AS condition by requiring its conclusion only over the outer inverse subfamily \emph{i.e.} for every outer inverse $x$ of $a$, some unit $u$ satisfies $a-u\in\mathrm{ann}_l(x)$ and see if stronger structural conclusions hold. %Together with the simple observation that left AS ring $\implies $ partially unit-regularity (where every von Neumann regular element is unit-regular), we have the following corollary. 
%As a consequence of the above theorem, we weaken the left (right) AS condition further: instead of assuming $Ra + \mathrm{ann}_l(b) = R$ for some $b$, we weaken to the case when $a$ admits an outer inverse (in which case there exists $x \in R$ such that $Ra + \mathrm{ann}_l(x) = R$) for arbitrary $a$ and show stronger structural conclusions hold.

\begin{cor}{\label{prp7}}
	For any ring $R$ the following conditions are equivalent :
	\begin{enumerate}[label=(\roman*)]
		\item $R$ is partially unit-regular. 
		\item For every $a \in R$ and any outer inverse $x \in T(a)$,   there exists an invertible element $u \in U(R)$ such that $ a- u \in \annr(x) \cap \annl(x)$
		\item For every  $a \in R$ and any outer inverse $x \in T(a)$,   there exists an invertible element $u \in U(R)$ such that $ a- u \in \operatorname{ann}_{l}(x) \left(\text{equivalently, } \operatorname{ann}_r(x)\right)$
	\end{enumerate}
\end{cor}
\begin{proof}
$(i) \implies (ii):$	Let $a\in R$ and  $x$ be an arbitrary outer inverse of $a$. Then $xax = x$. Then $x$ is von Neumann regular, with $a$ as an inner inverse, and hence $x$ is unit-regular by $(i)$. Moreover $a' = axa$ is an outer inverse of $x$, using $xax = x$, and therefore  $a' \in T(x)$. Applying Theorem~\ref{urr}$(ii)$ to the unit-regular element
$x$ and its outer inverse $a'$, there exists $s \in U(R)$ such that
\[  a' + (1-a'x)s(1-xa') = axa + (1-ax)s(1-xa) \in U(R). \]
	Label $u = axa +(1-ax)s(1-xa)$.  Since $x(1-ax) = 0$ and $(1-xa)x = 0$, we obtain
\[ xu = x(axa) = xa \qquad \text{and} \qquad ux = (axa)x = ax, \]
hence $x(a-u) = 0 = (a-u)x$, that is, $a-u \in \annr(x) \cap \annl(x)$.
%Then by Theorem \ref{outer0}, there exists a reflexive inverse $x' \in Ref(a)$ such that $x \leq^{-} x'$ \emph{i.e.} there exists $p \in R$ such that $xpx = x$ and $px = px', xp = x'p.$ \\
	%Now $x'\in Ref(a)$ being von Neumann regular, is unit-regular by our assumption and since $ a \in T(x'),$ we have by Theorem \ref{urr} $(ii)$:
	%\[a +(1-ax')s(1-x'a) \in U(R)\]
	
	%Label $u = a +(1-ax')s(1-x'a)$. Therefore, $(a-u) \in \annl(x') \cap \annr(x')$. \\

%\noindent Now $(a-u) \in \annl(x') \implies (a-u)x' = 0 \implies (a-u)x'p = 0  \implies (a-u)xp = 0.$  Post-multiplying by $x$ and using $xpx = x$ we have \[ (a-u)x = 0 \implies (a-u) \in \annl(x)\]
%Similarly, using $px' = px$ we get $(a-u) \in ann_r(x)$ and subsequently, \[(a-u) \in \annl(x) \cap \annr(x)\]
	$(ii) \implies (i):$ 
	Let $a \in R$ be von Neumann regular and choose $x \in \mathrm{Ref}(a)$. 
	In particular, $x$ is an outer inverse of $a$, so by hypothesis there 
	exists $u \in U(R)$ with
	$$a - u \in \mathrm{ann}_\ell(x) \cap \mathrm{ann}_r(x).$$
	From $x(a-u) = 0$ and $(a-u)x= 0$ we obtain $xa = xu$ and $ax = ux$. Since  $x \in \operatorname{Ref}(a)$ we must have   $a \leq^{-} u$ where $u \in U(R).$  \\
Consequently,  $a$ is unit-regular. \\
$(ii) \iff (iii):$ Straightforward. \\
$(iii) \iff (i):$ Follows directly from \cite{siddique1} Prop. $9(iii)$. \\
Equivalently, the proof is similar if we  replace $\operatorname{ann}_{l}(x)$ with $\operatorname{ann}_{r}(x).$
\end{proof}
\medskip 

\noindent The above corollary alongwith Lemma \ref{ASouter} has some interesting consequences. Firstly,  using the characterization of Canfell for AS rings namely $R$ is left UG $\iff$ $R$ is left AS [Canfell {\cite[Cor.~4.4]{canfell}}],  we obtain the following theorem  of Marks.

%\begin{theorem}[Canfell {\cite[Cor.~4.4]{canfell}}]\label{thm:canfell}
%	The following statements are equivalent for a ring $R$.
%	\begin{enumerate}[label=(\roman*)]
	%		\item If $Ra+ann_l(b)=R$ for $a,b\in R$, then there exists a unit $u\in U(R)$ such that $a-u\in ann_l(b)$.
	%		\item $R$ is left uniquely generated (left UG).
	%\item If $Ra=Rb$ for $a,b\in R$, then $a=vb$ for some left unit $v\in R$.
	%	\end{enumerate}
%\end{theorem}

\begin{cor}[Marks {\cite[Theorem 20]{marks}}]
	A regular ring $R$ is unit-regular if and only if it is left UG.
\end{cor}
\begin{proof}
	%	For the converse, we use Lemma \ref{ASouter} to show that if $R$ is unit regular then $\{b : Ra + \annl(b) = R\} = \bigcup_{x \in T(a)} xR$ and therefore from Corollary \ref{prp7}, $R$ is AS $\iff$ R is UG. 
	$\implies :$  Using  the universal  inclusion $\{b : Ra + \annl(b) = R\} \supseteq \bigcup_{x \in T(a)} xR$ and Corollary \ref{prp7} on a von Neumann regular ring $R$, we obtain $R$ is unit-regular. \\
	To prove the converse, suppose $R$ is unit regular. If possible let $Ra + \annl(b) = R.$ Then from Lemma \ref{ASouter}, $b = xr$ for some $ x \in T(a)$. Now from Corollary \ref{prp7}, if there exists $u \in U(R)$ such that $a-u \in \operatorname{ann}_{l}(x) \implies (a-u)x = 0 \implies (a-u)b = 0.$ Thus $R$ is left AS $\iff$ R is left UG ([Canfell {\cite[Cor.~4.4]{canfell}}])
\end{proof}

Secondly,  using the same characterization of Canfell, Lemma \ref{ASouter} and Corollary \ref{prp7},  we observe the following well-known fact  using  the left-right symmetric notions of exchange  and partial unit-regularity. %and the opposite ring $R^{\operatorname{op}}$.  %in an exchange ring $R$, left UG is equivalent to right UG [Khurana, Lam {\cite[Thm.~6.5]{lam4}}].
\begin{theorem}[Khurana, Lam {\cite[Thm.~6.5]{lam4}}]
An exchange ring $R$ has left UG iff it has right UG.
\end{theorem}
%\begin{proof}
%	The  proof is straightforward bcfut we show it here for the sake of completition. \\
%	Suppose  $R$ is left UG $\iff$ R is left AS by Canfell's characterization. If possible, let $aR + \operatorname{ann}_{r}(b) = R$ for some $b \in R$. Then, by Remark \ref{ASouter1}, $b \in Rx $ for some $x \in  \operatorname{T}(a)$. Then 
%\end{proof}

Thirdly, using $(ii)$ of our Corollary \ref{prp7} we  combine the two criteria in  Corollary 3.3  in Chen ~\cite{Chen Stable Range} and weaken  reflexive inverses with outer inverses. %where $a^{+} \in \operatorname{Ref}(a)$.
\begin{cor}
Let $R$ be an exchange ring. Then following conditions are equivalent:
	\begin{enumerate}[label = (\roman*)]
		\item $R$ has stable range one. 
		\item For any regular $a \in R$, and every $a^{+} \in T(a)$, there exists  $e\in \operatorname{ann}_{r}(a^{+}) \cap \operatorname{ann}_{l}(a^{+}) $ and $u\in U(R)$ such that
		\[
		a = e + u.
		\]
	%	\item For any regular $a \in R$, there exists  $e\in \operatorname{ann}_{l}(a^{+})$ and $u\in U(R)$ such that
	%	\[
	%	a = e + u.
	%	\]
	\end{enumerate}
\end{cor}
\begin{remark}
It can also be seen that the above  consequences is equivalent for every reflexive inverse $a^{+} \in Ref(a)$.
\end{remark}

The preceding result characterizes unit-regularity of an element by imposing unit-regularity on its outer inverses. If we weaken this requirement, such as instead of asking that the relevant outer inverse be unit-regular, we ask only that it  be an inner inverse of $a$, then we get the usual characterization of von Neumann regularity: \emph{an element $a\in R$ is regular if and only if there exists $x\in R$ such that
\[axa=a
\qquad\text{and}\qquad
xax=x.\]}
Infact, using Theorem \ref{outer} and  Lemma \ref{refinement} together with  Theorem \ref{partialorder} $(v)$, one obtains the following proposition.

\begin{prop}
	An element $a$ in a ring $R$ is von Neumann regular if and only if there exists outer inverse $y$ of $a$ and  reflexive inverse $x$ of $a$ such that 
%	\[x -y = (1-yy^{+})s(1-y^{+}y)\] for some $y^{+} =  \in  Ref(y)$.\\ 
	\[x -y = \left((1-ya) - r_{\scriptsize 1}(1-ya) \right) s \left((1-ay) - (1-ay)r_{\scriptsize 2} \right) \] 
	for some $r_{\scriptsize 1},r_{\scriptsize 2}, s \in R$\\
%	Additonally, $x-y$ is invertible if and only if  $a$ is unit-regular.\\
\end{prop}
%\begin{proof}
%	$(\Longleftarrow): Straightforward$ \\
%	$(\Longrightarrow):$ 	The first part of the proof is obtained by concatenating Lemma \ref{refinement}, Theorem \ref{outer}, and  \cite{hartwig} Proposition $3$. For completion, we show it below.\\
%	Choose an arbitrary outer inverse $y$ of $a$. Then by Theorem \ref{outer} there exists a reflexive inverse $x$ of $a$ such that $y \leq^{-} x$. Then by  \cite{hartwig} Proposition $3 (i)$ there exists some $y^{+} \in Ref(y)$ such that 
%	\[x -y = (1-yy^{+})s(1-y^{+}y)\] 
%	Using Lemma \ref{refinement}, we can substitute $y^{+} = \left[a + (1-ay)r\right](ya).w$ to obtain 
%	\[x -y = \left((1-ya) - r_{\scriptsize 1}(1-ya) \right) s \left((1-ay) - (1-ay)r_{\scriptsize 2} \right) \] 
 %	for some $r_{\scriptsize 1},r_{\scriptsize 2}, s \in R.$\\
 %	The converse is straightforward.
%\end{proof}
\noindent The term $(x-y)$ turns out to be an important factor in the characterization of broad classes of rings. In fact, continuing in this direction, we observe  an interesting analogue  when  $x,y$ are  assumed as outer inverses of $a, (1-a)$ respectively such that $ax + (1-a)y = 1$. %Knowing that whenever $a$ is von Neumann regular then both $a$ and $(1-a)$ are suitable (exchange),  Nicholson’s criterion gives idempotents in $aR$ and $(1-a)R$, and those representatives can be replaced by outer inverses while preserving the products $ ax$ and $(1-a)y.$  \\%since whenver an element $a$ is suitable, there exists idempotents $e \in Ra$ and $(1-e) \in R(1-a)$ as per Nicholson's criterion. We can subsequently choose $x, y \in R$ such that $ e = xa$ and $(1-e) = y(1-a)$ where $xax = x$ and $y(1-a)y= y.$ \\
 
\noindent  First we consider the following lemma. 
 \begin{lem}{\label{outerinverse}}
Let $a \in R$ be von Neumann regular with $z \in \operatorname{I}(a)$, and let $x$ be
an outer inverse of $a$. Label $B_x= \{\, y \in \operatorname{T}(1-a) \;:\; ax + (1-a)y = 1 \,\}$.

Then the following are equivalent:
 \begin{enumerate}[label = (\roman*)]
 	\item $(1-ax) \in (1-a)R$
 	\item  $B_x \neq \phi$.
 \end{enumerate} 
  In this case, $B_x = \{(1-az) + r(1-ax) : r \in \operatorname{ann}_{r}(1-a) \cap R(1-ax)\}$
 \end{lem}

\begin{proof}
	\noindent  Clearly for every outer inverses $x$ of $a$, $(1-ax)$ is an idempotent. \\
$(\Longrightarrow)$ Suppose $(1-ax)= (1-a)z$ for some $z\in R$. Choose $y =  z(1-ax).$ Then $y(1-ax) = y$ and $y(1-a)y = y(1-a)z(1-ax) = y(1-ax) = y.$ Thus $y \in B_x$. \\
The converse is straightforward.\\

%\noindent $(\Longleftarrow)$ If there exists an outer inverse $y$ of $(1-a)$  such that $y(1-ax) = y$
In addition, suppose  $B  \neq \phi$  and choose arbitrary  $ y \in B$. Label 
$y_{\scriptscriptstyle 0} = 1-az$. It satisfies $y_{\scriptscriptstyle 0}(1-a)y_{\scriptscriptstyle 0}=y_{\scriptscriptstyle 0}$ since $(1-az)a = a-aza=0$. Furthermore 
\[
(1-a)y_{\scriptscriptstyle 0} = (1-a)(1-az) = (1-ax)
\]
%and $y_0$ is idempotent.
Set $k = y - y_{\scriptscriptstyle 0}$. Then clearly, $(1-ax)k = (1-ax) (y-y_{\scriptscriptstyle 0}) = (1-ax) -(1-ax) = 0$ and $k (ax) =  (y-y_{\scriptscriptstyle 0})ax  = 0 \implies k = k(1-ax) \in R(1-ax)$. Thus \[B \subseteq \{(1-az) + k(1-ax) : k \in \operatorname{ann}_{r}(1-ax) \cap R(1-ax)\} \]
Conversely, for any element $y_{\scriptscriptstyle k} \in \{(1-az) + k(1-ax) : r \in \operatorname{ann}_{r}(1-a) \cap R(1-ax)\} $ we clearly have
\begin{align*}
(1-a)y_{\scriptscriptstyle k} &= (1-a)\underbrace{(1-az)}_{ y_{\scriptscriptstyle 0}}+ \underbrace{(1-a)k}_{0}(1-ax) \\
&= 1-ax
\end{align*} 
and $y_{\scriptscriptstyle k}(1-a)y_{\scriptscriptstyle k}= y_{\scriptscriptstyle k} (1-ax) = y_{\scriptscriptstyle k}$. Thus \\
\[B_x= \{(1-az) + k(1-ax) : k \in \operatorname{ann}_{r}(1-a) \cap R(1-ax)\}\] 
\end{proof}

\begin{prop}{\label{reflex1}}
		Let $R$ be a ring and $a \in R$.  Then $a$ is von Neumann regular  if and only if there exists a reflexive $($outer$)$  inverse $x$ of $a$ and outer inverse $y$ of $(1-a)$ such that the following hold.
	\begin{enumerate}[label = (\roman*)]
		\item $ax +(1-a)y = 1$
		\item $a(x-y)a = a$
	\end{enumerate}
	Additonally, if $(x-y)$ is invertible then $a$ is clean (and unit-regular).
\end{prop}

\begin{proof}
	$(\Longleftarrow): Straightforward$ \\
		$(\Longrightarrow):$ 
	Let $a$ be von Neumann regular. Choose  $z$ to be an inner-inverse of $a$. Then by the Refinement Lemma \ref{refinement}, the set of all reflexive inverses of $a$ is given by 
	\[\text{Ref}(a) = \{ (z+ (1-za)r)(az)(1+ azt(1-az)): \forall r,s,t \in R\}\]
	
	\noindent 	Choose a reflexive inverse of $a.$  Then, it is of the  form 
	\[x = \left(z + (1-za)r \right) \left(az\right)\left(1+azt(1-az)\right) \text{ for some } r,t \in R\]
	\noindent Then, $ax = az\left(az + azt(1-az)\right) = az + azt(1-az)$ for some $t\in R$. \\
	
\noindent  Furthermore, $1-ax = 1 - (az + azt(1-az)) = (1-azt)(1-az).$ Pluggin in $t = a$ gives us $1-ax = (1-a)(1-az)$. \\
\noindent 	We find $y \in R$ such that $y(1-a)y = y $ and $(1-a)y = 1 - ax.$ \\
\noindent Since $(1-ax) \in (1-a)R$, from Lemma \ref{outerinverse}, the set $B_x= \{y \in \operatorname{T}(1-a): (1-a)y + ax = 1 \}  \neq \phi,$ and therefore  we can get a family of such outer inverses  of $(1-a)$ given by \[ \{(1-az) + k(1-ax) : k \in \operatorname{ann}_{r}(1-a) \cap R(1-ax)\}\]
	
	\noindent Therefore, choosing $y_{\scriptscriptstyle k} =   (1-az) \, +  k(1-ax) \in B$  and $x_{\scriptscriptstyle r} =  x= (z + (1-za)r ) (az)(1+a(1-az))$  for any $r \in R$  satisfies 	\begin{enumerate}[label = (\roman*)]
		\item $x_{\scriptscriptstyle r}ax_{\scriptscriptstyle r} = x_{\scriptscriptstyle r}$  and $ax_{\scriptscriptstyle r}a= a$
		\item $y_{\scriptscriptstyle k}(1-a)y_{\scriptscriptstyle k}=  y_{\scriptscriptstyle k}$
		\item $ax_{\scriptscriptstyle r} +(1-a)y_{\scriptscriptstyle k}= 1$\\
	\end{enumerate}
	
To prove the additional part, we use the term $(x-y)$ to characterize $a$ explicitly. 
	\noindent Set $e = ax$. Then $e^2 = e, \; x \cdot e = x,  \text{ and }y \cdot(1-e) = y. $    Then we can write, 
	\begin{align*}
		a(x-y) - (x-y)(1-e) &= ax - ay - (x-y) + (x-y)e \\
		&= (ax) + ((1-a)y )- x + xe \\
		&= (e) + (1-e)\\
		&= 1
	\end{align*}
	This is true for all values of $x=  \left(z + (1-za)r \right) \left(az\right)\left(1+a(1-az)\right) \text{ for  } r \in R$. \\
	
	\noindent Now, if $(x-y) = u$ is invertible for any  value of $r \in R,$ then we can write $a = u(1-e)u^{-1} + u^{-1}$, which shows that $a$ is clean. Additionally,\\
	
	\noindent $axa = a \implies (1-e)a = 0 \implies y(1-e)a = 0 \implies y \cdot a = 0$. Consequently, $aua = a(x-y)a = axa = a.$
\end{proof}

Continuing in this direction, we next consider what remains of the preceding characterization when the inverse data are weakened. In Proposition \ref{reflex1}, regularity is obtained from a reflexive inverse of $a$, an outer inverse of $1-a$, and a compatibility condition involving their difference. If the reflexive inverse is replaced by an outer inverse and this compatibility condition is omitted, we find that the resulting condition is no longer regularity but suitability (exchange).

\begin{theorem}\label{suitable}
	Let $R$ be a ring and $a \in R$.  Then the following hold:
	\begin{enumerate}[label = (\roman*)]
		\item $a$ is suitable
		\item There exists outer inverses $x$ of $a$ and $y$ of $(1-a)$ such that $ax + (1-a)y = 1$
		\item There exists outer inverses $x$ of $a$ and $y$ of $(1-a)$ such that  
		\[a(x-y) - (x-y)(1-e) = 1\]
		\item $(1-a)$ is suitable.
	\end{enumerate}
\end{theorem}
\begin{proof}
$(i) \implies (ii):$ Let $a$ be suitable. Then there exists an  idempotent $e \in aR$ and $(1-e) \in (1-a)R$ as per Nicholson's criterion.\\ Set $x, y \in R$ such that $ e = ax$ and $(1-e) = (1-a)y,$ which subsequently gives us $ax +(1-a)y = e+ (1-e) =1$.  We can choose $x,y \in R$ such that $xax = x$ and $y(1-a)y= y,$ for otherwise  set $x' = xax$. Then $x'ax' = x'$ and $ax' = ax \implies (1-a)y +ax' = 1$. \\
$(ii) \implies (i):$ Choose $x,y \in R$ such that $xax = x, y(1-a)y = y$ and $ax + (1-a)y = 1$. Set $e = ax$ and $f = (1-a)y$. Then $f = 1-e$ where $e \in aR$ and $f = 1-e \in (1-a)R$.\\
$(i) \implies (iii):$ 	 Let $a$ be suitable. Then there exists an idempotent $e^2 = e \in E(R)$ such that $e = ar$ and $1-e = (1-a)s$. Choose $r,s$ such that $re = r$ and $s(1-e)=s$. \\
Then $rar= r, s(1-a)s = s$ and 
\begin{align*}
	a(r-s) - (r-s)(1-e) &=	a(r-s) + s(1-e) \\
	&= ar - as + s \\
	&= ar + (1-a)s\\ 
	&= e + (1-e) =1
\end{align*}
$(iii) \implies (i)$ Let there exists $x,y \in R $ such that \[a(x-y) -  (x-y)f = 1 \] for some $f = f^2 \in R$. Post-multiplying by $(1-f)$ and $f$ respectively, we get
\[a(x-y)(1-f) = (1-f)  \hspace{3mm} \text{ and } \hspace{3mm} (1-a)(y-x)f = f  \]
Thus $a$ is exchange.\\
Finally, $(i)\Longleftrightarrow(iv)$ follows from the symmetry of Nicholson’s criterion.
\end{proof}
\begin{remark}
An equivalent criterion $(x-y)a-(1-e)(x-y)=1$ for suitability of elements can be obtained by using idempotents in $Ra$ and $R(1-a)$. 
\end{remark}
As a consequence, we give a similar characterization of clean rings in terms of outer inverses. After obtaining this criterion, we found that our clean case resembles  a result of Mary \cite{mary} which is proved  using an interesting Bott-Duffin $(e,f)$-inverse approach in \cite[Cor. ~2.3]{mary}. 

\begin{cor}
	Let $R$ be a ring and $a \in R$.  Then $a$ is clean if and only if there exists outer inverses $x$ of $a$ and $y$ of $(1-a)$ such that $ax + (1-a)y = 1$ and  $(x-y)$ is invertible.
\end{cor}
\begin{proof}
	$(\Longleftarrow)$	If there exists outer inverse $x$ of $a$ and $y$ of $(1-a)$ such that  $ax + (1-a)y = 1$ then by Theorem \ref{suitable} we can write 
		\[a(x-y) - (x-y)(1-e) = 1\] 
		Since $(x-y)$ is invertible, setting $u = x-y$  we can write $a = u(1-e)u^{-1} + u^{-1}$, which shows that $a$ is clean. \\
		
			$\Longrightarrow$ Let $a$ be clean. Then there exists an idempotent $f^2 = f \in R$ and unit $u \in U(R)$ such that
			\[a = f + u^{-1}\]
			Set $(1-e) = u^{-1}fu$. Then, $ f = u(1-e)u^{-1}$ and we can write \[a = u(1-e)u^{-1} + u^{-1} \hspace{5mm} \cdots \hspace{5mm} eq.(i)\]
			Set $ x = ue$ and $y = (-u)(1-e)$. Then, we can see that $x - y = u \in U(R)$.\\
			Infact, $xax = ue\left(u(1-e)u^{-1} + u^{-1}\right) ue = ue = x$. \\
			Now rerranging terms in $eq.(i)$ we  see that
			\begin{align*}
				1-a &= 1- \left(u(1-e)u^{-1} + u^{-1} \right) \\
				&= ueu^{-1} - u^{-1}
			\end{align*}
			Then, $y(1-a)y = (-u)(1-e)\left(ueu^{-1} - u^{-1}\right)(-u)(1-e)= (-u)(1-e) =y$
			\begin{align*}
				\text{Furthermore, } ax + (1-a)y &= \left(u(1-e)u^{-1} + u^{-1}\right)(ue) + \left(ueu^{-1} - u^{-1}\right)(-u)(1-e)\\
				&= e + (1-e) = 1
			\end{align*}
			Hence, there exists outer inverses $x$ of $a$ and $y$ of $(1-a)$ such that $(x-y)$ is invertible and  $ax + (1-a)y = 1$.\\
			Additionally, from $eq.(i)$ we can also write \[a (x-y) - (x-y)(1-e) = 1\]
\end{proof}

%In \cite{khurana2} theorem they prove solvability and show that the coefficient $x$ is an arbitrary element of $R$ in the suitable case and a unit in the clean case. In our formulation we find that the coefficient $x$ is not merely shown to exist, but can be realized as  the difference of  outer inverses  of $a$ and $(1-a)$ respectively.
\noindent It is useful to observe that the equations appearing in our characterizations have the form of an inhomogeneous Sylvester-type equation.  Indeed, the identity
$(x-y)a-(1-e)(x-y)=1$ may be viewed respectively as equation of the form
\[AX-XB=C\]
with $X=x-y, B = F^2 =F$ and $C=1$. 
and is well-known as the classical  Sylvester equation.   The Sylvester equation $AX-XB=C$, originally arising in matrix theory, has a substantial literature concerning its solvability, uniqueness, and block-similarity interpretations; For a detailed discussion of this equation, authors  refer to the exposition article of Bhatia-Rosenthal \cite{bhatia}.  In the ring-theoretic setting, equations of this same left-right linear form have also been used to study exchange phenomena.  Khurana, Lam, and Nielsen \cite{lam5} studied the equation
\[XA-BX=I\]
over a noncommutative ring and showed that, when $B$ is idempotent, its solvability is closely tied to Nicholson’s exchange, or suitable, elements.  In particular, they called the  equation \[xa-fx=1,\qquad f^{2}=f\]  \emph{exchange equation} for $a$, and its solvability characterizes suitable elements.  With additional restrictions on the solution $x$, the same framework also recovers characterizations of clean and strongly clean elements.

\begin{theorem}[\cite{lam5} Theorem 2.1]
For any element $a\in R$, the following hold.
\begin{enumerate}[label = (\roman*)]
	\item $a$ is suitable iﬀ there exist $f \in idem (R)$ and $x\in R$ such that $xa-fx= 1$.
	\item $a$ is clean  iﬀ there exist $f \in idem(R)$ and $x\in U(R)$ such that $xa-fx= 1$.
	\item $a$ is scn(R) iﬀ there exist $f\in idem (R)$ and $x \in CR(a)$ such that $xa -fx =1$.
\end{enumerate}
\end{theorem}

\noindent What is special in our present setting is the additional information imposed on the solution. We show that in each of these cases the equations are reversible and the solution $X$ is not merely the solvability of the linear equation of the form $AX-XB=C$; rather, it is explicitly realized as $X=x-y$, where $x$ and $y$ are outer inverses explicitly associated  with $a$ and $1-a$ respectively.  Hence its not just the appearance of the Sylvester-type equation alone, but the identification of its solution as a difference of prescribed pseudo inverses, together with the additional hypotheses required in the corresponding characterizations of regular, suitable, and clean elements.\\

 %In our setting $X$ is not an arbitrary unknown but is prescribed as the difference of two outer inverses,
\noindent We proceed in this fashion to obtain equivalent characterizations for \emph{unit-exchange} and \emph{unit-clean} elements as defined in \cite{siddique1}.  Unit-exchange may be viewed as the natural analogue of suitability (exchange) in which the identity  occurring in the usual exchange decomposition is allowed to be replaced by an arbitrary unit.\\
An element $a\in R$ is \emph{left unit-exchange} if there exist $e^2=e\in R$ and $u\in U(R)$ such that
\[
Re \subseteq Ra \qquad\text{and}\qquad R(1-e)\subseteq R(a-u).
\]
Dually, $a$ is \emph{right unit-exchange} if $fR\subseteq aR$ and $(1-f)R\subseteq (a-v)R$ for some idempotent $f$ and unit $v$. The notions of left- and right-unit exchange ring is symmetric.
% An element $a\in R$ is called \emph{unit-exchange} if there exist a unit $u\in U(R)$ and an idempotent $e\in R$ such that
%\[Re\subseteq Ra
%\qquad\text{and}\qquad
%R(1-e)\subseteq R(u-a).\]

Similarly, an element $a\in R$ is called \emph{unit-clean} if it can be written as
$a=eu+v$,
where $e^{2}=e$ and $u,v\in U(R)$.  Equivalently, $a$ is a sum of a unit-regular element and a unit.  %This places unit-clean elements in direct parallel with clean elements: the idempotent part is replaced by a unit-regular component, while the remaining summand is still required to be a unit.
\begin{prop}
	Let $R$ be a ring. Then the following are equivalent 
	\begin{enumerate}[label = (\roman*)]
		\item $a$ is unit-exchange.
		\item There exists  outer inverses $r,s$ of $ua$ and $1-ua$ respectively, such that 
		\[(r-s)ua- f(r-s)= 1\]
		for some idempotent $f \in R$ 
				\item For some $u \in U(R)$, there exists  outer inverses $r,s$ of $au$ and $1-au$ respectively, such that 
		\[au(r-s) - (r-s)f = 1\]
		for some idempotent $f \in R$ 
		\item $(1-a)$ is unit-exchange
	\end{enumerate}	
	Additionally $a$ is unit clean if and only if $(r-s)$ is invertible.
\end{prop}

\begin{proof}
	
	$(i) \implies (ii)$ The proof follows similarly as in the case of suitable rings above in Theorem \ref{suitable}. \\ %Let $a$ be unit-exchange. Then there exists a unit $u \in U(R)$ and idempotent $e^2 = e \in E(R)$ such that $e = aur$ and $1-e = (1-au)s$. Choose $r,s$ such that $re = r$ and $s(1-e)=s$. \\
	%Then $raur= r, s(1-au)s = s$ and 
	%\begin{align*}
	%	au(r-s) + (r-s)(1-e) &=	au(r-s) + s(1-e) \\
	%%	&= aur - aus + s \\
	%	&= aur + (1-au)s\\ 
	%	&= e + (1-e)\\
	%	&=1
	%\end{align*}
	$(ii) \implies (i)$ Conversely, let there exists $r,s \in R $ such that $au(r-s) -  (r-s)f = 1 $ for some $f = f^2 \in R$ and for some $u\in U(R).$  \\
	
Multiplying this identity on the left by $1-f$, we obtain
\[(1-f)(r-s)ua=1-f.\]
Define  $\gamma=(1-f)(r-s)u.$
Then $\gamma a=1-f \implies 1-\gamma a=f.$ Moreover, $\gamma$  is an outer inverse of $a$, since
\begin{align*}
\gamma a\gamma
&=(1-f)(r-s)u\,a(1-f)(r-s)u \\
&=(1-f)(r-s)u(1-f)=\gamma.
\end{align*}

\noindent Next, multiplying
$(r-s)ua-f(r-s)=1$
on the left by $f$, we get
$f(r-s)ua-f(r-s)=f.$
Equivalently,
$f(r-s)(ua-1)=f.$
Since
$ua-1=u(a-u^{-1}),$
this gives
$f(r-s)u(a-u^{-1})=f.$
Multiplying on the right by $f$, and using $f^{2}=f$, yields
\[f(r-s)u(a-u^{-1})f=f.\]
Thus, if
$\sigma=f(r-s)u
\qquad\text{and}\qquad
v=u^{-1},$
then
\[\sigma(a-v)(1-\gamma a)=1-\gamma a.\]

%	Post-multiplying by $(1-f)$ and $f$ respectively, we get 
%	\[au(r-s)(1-f) = (1-f)  \hspace{3mm} \text{ and } \hspace{3mm} (1-au)(s-r)f = f  \]
\noindent 	Then by \cite{siddique1}  Proposition 20, $a$ is unit-exchange. \\
For the additional part, the solution follows similarly as above in the case of suitable rings.
\end{proof}

%\begin{remark}
%Let $a$ be unit-exchange. Then there exists a unit $u \in U(R)$ and idempotent $e^2 = e \in E(R)$ such that $e = aur$ and $1-e = (1-au)s$. Choose $r,s$ such that $re = r$ and $s(1-e)=s$. 	If $r-s$ is invertible, then the outer inverses $r,s$ are unit-regular for 
%	\[r = (r-s)e  \hspace{3mm} \text{ and } \hspace{3mm} s = (r-s)(1-e) \]
	
%	\noindent In addition,  $a$ is unit-clean since by Theorem 35, there exists a unit $v = r-s$ and idempotent  $e = 1-f$ such that 
%	\[a(uv)e = e  \hspace{3mm} \text{ and } \hspace{3mm} (a-u^{-1})(uv)(1-e)= 1-e \] 
%\end{remark}

\bigskip

We now proceed to the case of strongly clean rings. An element $a\in R$ is called \emph{strongly clean} if it can be written as
$a=e+u$,
where $e^2=e, u\in U(R)$, and $e$ commutes with $u$ or $e \in C(u)$.\\
We first consider the following consequence
\begin{prop}{\label{sclean}}
	Let $a \in R.$ If there exists  outer inverses $x$ of $a$ and $y$ of $(1-a),$ such that $(1-a)y + ax = 1$ and  $(x-y) \in C(a),$ then $a$ is strongly clean.
\end{prop}

\begin{proof}
	Let $a \in R$. Choose $x,y \in R$ such that \[xax = x, \; y(1-a)y = y, \; \text{ and } ax + (1-a)y = 1\] Set $e = ax.$ Then from $ x(1-e)=0, \; y(1-e) = y, \; (1-e) = (1-a)y$  we have 
	\[a(x-y) - (x-y)(1-e) = 1\]
	Now, $x-y \in C(a) \implies  (x-y)\left(a - (1-e)\right) = 1$. \\
	On the other hand, 
	\begin{align*}
		\left(a - (1-e)\right) (x-y) &= a(x-y) - (1-e)x + (1-e)y \\
		&= \underbrace{ax +(1-a)y}_{1}  - ey - (1-e)x \\
		&= 1 - ey - (1-e)x  \hspace{35mm} \cdots \hspace{5mm}(i)\\
		%&= 1 + e(x-y)(1-e) - (1-e)(x-y)e  \hspace{5mm} \cdots \hspace{5mm}(ii)\\
	\end{align*}
	Since $(x-y)$ commutes with $a$ we have $(x-y)a  = a(x-y)$ where 
	\begin{align*}
		a(x-y) &= ax - ay \\
		&= e + (1-a)y - y \\
		&= 1 + (1-e) - 1 \\
		&= 1-y
	\end{align*}
	Therefore, $(x-y)a  = 1-y \implies (x-y)ax  = (1-y)x \implies yax = yx$ where $yax = y\cdot e = (y(1-e))\cdot e=0$. Therefore $yx = 0$. \\
	Then, $(1-e) = (1-a)y \implies  (1-e)x =(1-a)yx = 0 \implies x = ex$. Thus, $x = ex = xe \in eRe$.\\
	Then from \emph{eq.}$(i)$ above, since  $e(x-y)(1-e)$ is a nilpotent element we have
	\begin{align*}
		\left(a - (1-e)\right) (x-y) &= 1 - ey - (1-e)x  \\
		&= 1 + e(x-y)(1-e)  \in U(R)
	\end{align*}
	Thus, $(x-y)$ is both left and right-invertible and is therefore an unit in $R$ with \[(x-y)(a - (1-e) )= 1\]
	 Therefore  $a -(1-e) = (x-y)^{-1}$  where $(x-y)$ commutes with $a$ and hence, also with $e=ax.$ Consequently, $a$ is strongly clean. \\
	\begin{remark}
		For any $a \in R$ if there exists  outer inverses $x$ of $a$ and $y$ of $(1-a),$ such that $(1-a)y + ax = 1$ then from above we see that  $(x-y)$ commutes with $a \implies (x-y)$ commutes with $ax$. \\
		On the other hand, \textbf{if  $R$ is directly finite}, then  $(x-y)$ commutes with $ax$ if and only if $(x-y)$ commutes with $a$. This is because we already have \[a(x-y) - (x-y)(1-e) = 1 \] for $e = ax$. If  $(x-y)$ commute with $ax$ and $R$ is directly finite, then $[a - (1-e)](x-y) = 1 \implies (x-y)[a-(1-e)] =1 \implies (x-y)a = a(x-y)$.
		Therefore, 
		\[(x-y) \in C(a) \Longleftrightarrow (x-y) \in C(ax)\]
		The direction is not reversible in general. 	Let  $V$ be a $k$-vector space with  basis
		$\{v_0, v_1, v_2, \dots\}$. Set
		\[
		R \;=\; \operatorname{End}_k(V).
		\]
		This ring $R$ is not directly finite. %as the shift maps defined below satisfy
		%$TS = 1_V$ while $ST \neq 1_V$.  
		 Define the \emph{forward} and
		\emph{backward} shifts on $V$ by
		\[
		S : v_i \longmapsto v_{i+1} \ (i \geq 0),
		\qquad
		T : v_i \longmapsto v_{i-1} \ (i \geq 1), \quad T v_0 = 0,
		\]
		so that $TS = 1_V$ and $ST = 1_V - p_0$, where $p_0(v_0)=v_0,p_0(v_i)=0$ for $i\geq 1.$  \\
		
	\noindent 	Define $A,X,Y\in R$ given by
		\[
		A(v_0)=0= A(v_1)\qquad A(v_i)=v_{i-1}\quad (i\geq 2),
		\]
		\[
		X(v_0)=0,\qquad \qquad \hspace{6mm} X(v_i)=v_{i+1}\quad (i\geq 1),
		\]
		\[
		Y(v_0)=v_0,\qquad \qquad \hspace{5mm} Y(v_i)=0\quad (i\geq 1). \hspace{4mm}
		\]

		Then $ AX \;=\;  1_V  - p_0$. It can  easily be seen that $AX \in \operatorname{Idem}(R)$. Label $ E = AX$.\\
		Label $w = X-Y$ where
\[
w(v_0)=-v_0,\qquad w(v_i)=v_{i+1}\quad (i\geq 1).
\]
		Then the following identities are immediate:
%\[e := ax = TS \oplus 0 = 1_V \oplus 0 \quad(\text{projection onto } V \text{ along } ku), \qquad w := x - y = S \oplus (-1)\]
		\begin{enumerate}[label=(\roman*)]
			\item $AXA=A, \,\,\, XAX=X$ 
			\item $Y(1-A)Y = Y$ 
			\item $AX +(1-A)Y = 1$
		 \item Furthermore, $(X-Y)$ commutes with $AX$ since on each $v_i : i \geq 1 $
		\[(X-Y)(AX)(v_i)=(X-Y)(v_i)=X(v_i)-Y(v_i)=v_{i+1}-0=v_{i+1}, \]  
		\[(AX)(X-Y)(v_i)=(AX)(v_{i+1})=v_{i+1}, \] and on $v_0:$  $(X-Y)(AX)(v_0) =  0 = (AX)(X-Y)(v_0)$
	\end{enumerate}
	However, $(X-Y)$ does not commute with $A$  on $v_1 $. %and subsequently, $(aw)(v_0) = a(v_1) = v_0$, whereas
	%$(wa)(v_0) = w(0) = 0$.\\
%		In each of these  cases, $A$ is strongly clean. \\
	\end{remark}
%\begin{remark}
\end{proof}

\begin{theorem}
	A  ring $R$ is strongly-clean if and only if  there exists outer inverses $x$ of $a$ and $y$ of $(1-a),$ such that $(1-a)y + ax = 1$
	and $(x-y)$ commutes with $a.$  %or $ (x-y)$ commutes with $ax$(or $xa$)
\end{theorem}

\begin{proof}
	Let $R$ be strongly clean. Then there exists $v \in U(R)$ and idempotent $f^2 = f$ such that $ a = f + v$ and $fv = vf$. Now we will rewrite it in the form of Sylvester's equation $AX - XB = C$ for some $B,C$  \\Set $ 1-e = v^{-1}(f)v$ and $u= v^{-1} \in U(R)$. Then we have $a = u(1-e)u^{-1} + u^{-1}$ such that $ue = eu$. Now  we re-write it as \[au-u(1-e) = 1 \text { where } ue= eu  \hspace{3mm} \cdots \hspace{3mm} (i)\]
	Set $x = ue$ and $y = -u(1-e)$. Then $(x-y)=u \in U(R).$ We can also rewrite $eq.(i)$ as 
	\[au  = 1 + u(1-e) \hspace{5mm} \text{ and } \hspace{5mm} (a-1)u = 1 - ue \] 
	With the chosen values of $x$ and $y$ we have 
	\begin{align*}
		(i) \; \; xax & = ue\underbrace{(a)u}e \hspace{100mm}\\
		&= ue(1+ u(1-e))e \\
		&= ue = x
	\end{align*}
	\begin{align*}
		(ii) \; \; y(1-a)y & = (-u)(1-e)\underbrace{(1-a)(-u)}(1-e) \hspace{100mm}\\
		&= (-u)(1-e)(1-ue)(1-e)\\
		&= (-u)(1-e) = y
	\end{align*}
	\begin{align*}
		(iii) \; \; ax +(1-a)y & = aue + (1-a)(-u)(1-e) \hspace{100mm}\\
		&= aue + (1-ue)(1-e)\\
		&= e + (1-e) = 1
	\end{align*}
	Additionally, since $ue = (x-y)e = e(x-y)$ where $u= (x-y) \in U(R)$ we see that 
	\[\left(a- (1-e)\right)(x-y) = 1 \; \;  \Longleftrightarrow \; \; (x-y)\left(a-(1-e)\right)= 1\]
	and consequently, $(x-y)$ commutes with $e = aue = ax$  $\Longleftrightarrow (x-y)$  commutes $a$.\\
	The converse follows from Proposition \ref{sclean}.
	
%	Conversely, for any element $a \in R,$ let there exists $x,y$ such that the above $(4)$ conditions are satisfied (in particular $(x-y)$ commutes with $a$ (or $ax$ or $xa$). We prove that $a$ is strongly-clean. \\
%	Clearly, $yax = 0$ and $ax + (1-a)y =1$. If $(x-y)$ commutes with $a$ then by Theorem , $a$ is strongly clean.\\
	%On the other hand suppose for any element $a \in R,$ let there exists $x,y$ such that the above $(4)$ conditions are satisfied and $(x-y)$ commutes with $ax.$  Set $e = ax$. Then by our assumption $(x-y)$ commutes with $e,$  and $x\cdot e = x, \; y(1-e) = y$.\\
	
	%\noindent We then have a sylvester's equation for the idempotent $B = (1-e), C = 1$ and the solution $X = (x-y)$ as follows : 
	%\begin{align*}
	%	a(x-y) - (x-y)(1-e)
	%	&= ax - ay - (x-y) + (x-y)e \\
	%	&= \bigl(ax + (1-a)y\bigr) - (x-xe) - ye \\
	%	&= 1 - (x-xe) - ye = 1  \qquad \cdots \qquad \because\ ye=0=x(1-e) \\
	%	\implies\ \bigl(a-(1-e)\bigr)(x-y) &= 1 \qquad\qquad\qquad\qquad\qquad \cdots \qquad \because (x-y)e = e(x-y)
	%\end{align*}
	
\end{proof}

As a consequence we recover the following corollaries that is given in Mary \cite[Cor. ~5.6]{mary} and Drazin \cite[Thm. ~5.5]{drazin} . 
\begin{cor}
	An element $a \in R$ is strongly clean if and only if there exists outer inverses $x$ of $a$ and $y$ of $(1-a)$   such that $ax +(1-a)y = 1 \,\,( \text{or } xa + y(1-a)=1)$ and $x,y \in C(a).$
\end{cor}
\begin{proof}
 If $x,y \in C(a)$ then $(x-y) \in C(a)$ and the proof follows from Theorem 37.\\
 On the other hand, if there exists outer inverses $x$ of $a$ and $y$ of $(1-a)$ such that $(x-y)\in C(a)$ then $(x-y)$ is invertible, and $a$ is strongly clean where  \[a-(1-e) = (x-y)^{-1}  \]
  Now from the earlier remark, $(x-y)$ also  commutes with $e=ax$. Therefore, $xy = yx = 0 \implies x \in C(x-y) \implies x \in C(\left(x-y\right)^{-1})$.\\
 Then using $a-(1-e) = (x-y)^{-1}$  we observe that  \[ \implies x\left(a-(1-e)\right) = \left(a-(1-e)\right)x \implies ax = xa \]
Similarly, we also obtain $y \in C(a)$.\\

\end{proof}

\section*{Acknowledgement}
\noindent This work was supported in part by a sabbatical leave granted by the University of Wisconsin - Eau Claire.

\end{document}